# A new class of high-order summation by parts finite-difference schemes

**L. Dovgilovich[1], R. Maksyutov[2], I. Sofronov[1,2]**

[1] Schlumberger, Moscow, Pudovkina 13

[2] MIPT, Moscow region, Dolgoprudny, Institutskii per.

ldovgilovich@slb.com , isofronov@slb.com

**Abstract.** We develop summation by parts (SBP) approach for generating high-order finite-difference schemes on the interval and propose new sets of schemes up to the 12[th] order. The coefficients of the schemes are governed by values of grid spacing near the ends of the interval: we shift one or two or three nodes at the ends of originally equidistant grid. The new finite-difference operators use forward and backward differences to avoid saw-tooth spurious solutions. Test computations point out two schemes (of 8[th] and 10[th] orders) that have the best accuracy among others; we yield the coefficients of new schemes.

## Introduction

Using the finite difference schemes based on SBP approach can produce efficient computational algorithms for the simulation of wave processes on rectangular grids [1, 3-8]. A known issue of SBP approximations of the first and second derivatives with a diagonal matrix of normalization operator is twofold reduction in the accuracy order of the difference scheme in the near boundary nodes (NBNs) [1, 5-8]. In addition, the maximum approximation order does not exceed $2p = 8$ for known schemes [1, 5-8] so far. An increase of the accuracy order can be achieved, for example, by the bigger number of NBNs with half order, i.e. $p$. For instance, the 10th order has been achieved in [4]; but it led to strong growth of the eigenvalue maximum modulus of the matrix differential operator, and, therefore, strong decrease of the Courant number of the correspondent difference scheme for solution of wave problems. Also, as proposed in [6, 7], an increase of the maximal order to *2p=10, 12, 14* in internal grid nodes (IGNs) can be achieved by keeping much lower order, the *fourth*, in *2p* NBNs.

In this paper, we use another approach to improve the accuracy of SBP schemes. Instead of the uniform grid we consider a grid with shifted NBNs: one, two, and up to three nodes at the boundary. By doing so we aim also avoiding the significant increase of the eigenvalue maximum modulus of the difference operator. As the result, we do improve the accuracy at the boundary for $8^{th}$ order schemes, and propose novel $10^{th}$ and $12^{th}$ order approximation schemes in IGNs. The SBP scheme [6-8] based on forward and backward differences for approximating the first derivative is used as the original SBP scheme for further optimization.

It should also be noted that because of the rather complex optimization algorithm for finding the grid nodes and coefficients of the proposed difference schemes, it was not possible to realize the full potential of the opportunities of the considered approach for high-order approximation (10th, 12th) in this study. But even these first schemes have quite high accuracy, not available on the uniform grid.

## 1 Difference schemes

### 1.1 Problem formulation

Consider the problem of finite-difference approximation of the second derivative operator $d^2/dx^2$ on the interval $[0,1]$.

We will use $p = 2, 3, ...$ for notation of the approximation order $2p$. Introduce the grid

$$\{x_0 = 0,\ x_1 = h_1,\ x_2 = x_1 + h_2, ..., x_i = x_{i-1} + h_i,\ x_N = 1\} \tag{1}$$

such that

$$h_p = h_{p+1} = ... = h_{N-p+1} = h,\ h_i = h_{N-i+1}, i = 1, ..., p-1,\ \sum_{i=1}^{N} h_i = 1,\ h_i > 0 \tag{2}$$

where $h$ is the constant spacing for IGNs. Thus the grid is symmetric, equidistant inside, and possibly non-equidistant for $p-1$ steps near the ends, see Fig. 1.

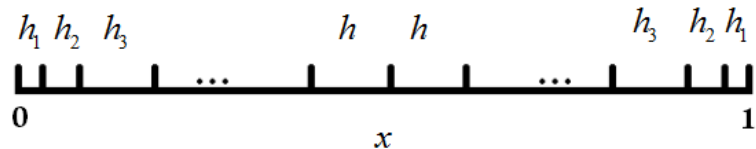

**Fig. 1 Computational grid**

## 1.2 SBP approach

Let us have difference operators $D^+$ and $D^-$ approximating the first derivative $d/dx$ on (1) in all points including boundary ones. We describe their construction hereafter. We introduce the scalar product in the space $\mathbb{R}^{N+1}$ of grid functions $\bar{u} = \{u_i\} = \{u(x_i)\}$, $i = 0...N$, on (1):

$$(\bar{u}, \bar{v})_H = \bar{u}^T H \bar{v}, \tag{3}$$

where $H > 0$ is a positive definite diagonal matrix (the entries are also defined hereafter). We require that the grid operators satisfy the Summation by Parts condition (SBP):

$$(\bar{u}, D^+\bar{v})_H + (D^-\bar{u}, \bar{v})_H = u_N v_N - u_0 v_0, \quad \forall \bar{u}, \bar{v} \in \mathbb{R}^{N+1}. \tag{4}$$

We will also use the matrix form of (4)

$$(D^+)^T H + H D^- = Q \tag{5}$$

where $Q = diag(-1, 0, \ldots 0, 1)$ is the diagonal $(N+1) \times (N+1)$ matrix.

**Equidistant grid**

First, let us consider the case $h_1 = \ldots = h_{p-1} = h$ in (1). The following results are known.

**A)** The operators $D^+ = D^- = D$ and $H$ are derived in [1] for values $2p = 2, 4, 6, 8$ such that

**A1)** the operator $D$ is the conventional first derivative *central difference* in IGNs $i = 2p, \ldots, N - 2p$ of $2p$ approximation order on the stencil $S_i : j = i - p, \ldots, i + p$

**A2)** the order of $D$ drops to $p$ in the remaining $2p$ nodes at the left and right ends.

**B)** The operators $D^+$, $D^-$ and $H$ are derived in [6-8] for values $2p = 2, 4, 6, 8$ such that

**B1)** the operators $D^+$, $D^-$ are conventional first derivative *forward* and *backward differences* in IGNs $i = 2p, \ldots, N - 2p$ of $2p$ approximation order on the stencils $S_i^+ : j = i - p + 1, \ldots, i + p + 1$ and $S_j^- : j = i - p - 1, \ldots, i + p - 1$, respectively;

**B2)** the order of $D^+$, $D^-$ drops to $p$ in the remaining $2p$ nodes at the left and right ends.

**C)** The operator $D_2$ is derived in [5] for immediate approximation of $d^2/dx^2$ at the values $2p = 2, 4, 6, 8$ such that

**C1)** the operator $D_2$ is the conventional second derivative *central difference* in IGNs $i = 2p, ..., N - 2p$ of $2p$ approximation order on the stencil $S_i : j = i - p, ..., i + p$

**C2)** the order of $D_2$ drops to $p$ in the remaining $2p$ nodes at the left and right ends.

Remark. The condition (4) from [8] is a generalization of conventional SBP [1] and permits us using forward and backward stencils in $D^-D^+$ to damp spurious saw-tooth oscillations in solutions.

**Grid with shifted NBNs**

Let us take an integer $p_g$, $1 \le p_g \le p$, and use a non-equidistant grid in $p_g - 1$ intervals of NBNs, i.e., we introduce free parameters $h_1, ..., h_{p_g - 1}$ while finding matrix coefficients of $D^+$ and $H$ ($D^-$ is defined from (5)).

Here IGNs are the nodes with $i = 2p, ..., N - 2p$.

Formulate the following problem.

**P.** Generate operators $D^+$, $D^-$ and $H$ for $2p \ge 4$ such that

**P1)** the operators $D^+$, $D^-$ are conventional first derivative *forward* and *backward differences* in IGNs of $2p$ approximation order on the stencils $S_i^+ : j = i - p + 1, ..., i + p + 1$ and $S_j^- : j = i - p - 1, ..., i + p - 1$, respectively (like in **B**);

**P2)** the order of $D^+$, $D^-$ is not less than $p$ in the remaining $2p$ nodes at the left and right ends; approximation error of the second derivative by $D^-D^+$ is *small as possible* for polynomials up to the $2p$ order;

**P3)** the maximal eigenvalue of the problem $(D^+)^T HD^+ \bar{u} = \lambda^2 H\bar{u}$ is *not too large* comparing to the maximal eigenvalue of the similar problem for *submatrices* of $(D^+)^T HD^+$ and $H$ in IGNs on the subspace of periodic functions.

Let us comment the conditions of the above problem. Firstly, we want to get the operators of a higher order than $2p = 8$ (known to date). Secondly, we would like to significantly improve the accuracy of approximation in NBNs (as we will see below, the higher order than $p$ is not possible). Thirdly, we would like to avoid a strong drop of the Courant number in NBNs while using the operator $D^-D^+$ to solve wave problems (i.e., comparing to IGNs where the operators are defined according to **P1**). Connection of $(D^+)^T H D^+$ with the operator of the second derivative $D^-D^+$ follows from (5): $D^-D^+ = -H^{-1}(D^+)^T H D^+ + H^{-1} Q D^+$. Therefore, for example, for a homogeneous Neumann problem the term $H^{-1} Q D^+$ vanishes, and for homogeneous Dirichlet problem it is simply not used as the endpoints are discarded from consideration.

### 1.3 Generation of coefficients for the matrix operators $D^+, D^-, H$

Here we describe implementation of **P1** and first part of **P2**.

Let us consider structure of matrices of our operators. Denote by $[A]_{M \times N}$ a matrix $A$ with $M$ rows and $N$ columns.

The matrix $H$ has the form:
$$[H]_{(N+1) \times (N+1)} = diag(\mu_1, \ldots, \mu_{2p}, 1, \ldots, 1, \mu_{2p}, \ldots, \mu_1) \times h \ , \tag{6}$$

where $\mu_i > 0$.

The matrix $D^+$:
$$[D^+]_{(N+1) \times (N+1)} = \frac{1}{h} \begin{pmatrix} D_l^+ & D_{lc}^+ & [0]_{2p \times (N+1-4p)} \\ [0]_{(N+1-4p) \times (p+1)} & D_c^+ & [0]_{(N+1-4p) \times (p-1)} \\ [0]_{2p \times (N+1-4p)} & D_{rc}^+ & D_r^+ \end{pmatrix}, \tag{7}$$

where

$$[D_l^+]_{2p \times 2p} = \begin{pmatrix} d_{1,1} & \cdots & d_{1,2p} \\ \vdots & \ddots & \vdots \\ d_{2p,1} & \cdots & d_{2p,2p} \end{pmatrix}, \tag{8}$$

$$[D_c^+]_{(N+1-4p) \times (N+1-2p)} = \begin{pmatrix} \beta_{-p+1} & \cdots & \beta_0 & \cdots & \beta_{p+1} & 0 & \cdots & 0 \\ 0 & \beta_{-p+1} & \cdots & \beta_0 & \cdots & \beta_{p+1} & \cdots & 0 \\ & & & & \ddots & & & \\ 0 & \cdots & 0 & \beta_{-p+1} & \cdots & \beta_0 & \cdots & \beta_{p+1} \end{pmatrix}, \tag{9}$$

$$[D^+_{lc}]_{2p \times 2p} = \begin{pmatrix} 0 & \cdots & & & & 0 \\ & \ddots & & & & \\ 0 & \cdots & & & & \vdots \\ \alpha_{p,1} & 0 & \cdots & & & \\ \vdots & \ddots & 0 & & & \\ \alpha_{2p,1} & \cdots & \alpha_{2p,p+1} & 0 & \cdots & 0 \end{pmatrix}, \qquad (10)$$

$$[D^+_{rc}]_{2p \times 2p} = \begin{pmatrix} 0 & \cdots & 0 & \gamma_{1,p-1} & \cdots & \gamma_{1,1} \\ \vdots & & & 0 & \ddots & \vdots \\ & & \ddots & & 0 & \gamma_{p-1,1} \\ & \cdots & & & & 0 \\ & \ddots & & & & \\ 0 & \cdots & & & & 0 \end{pmatrix}, \qquad (11)$$

and $[D^+_r]_{2p \times 2p}$ consists also of $2p$ stencils, each with $2p$ coefficients.

The similar form has the matrix

$$[D^-]_{(N+1) \times (N+1)} = \frac{1}{h} \begin{pmatrix} D^-_l & D^-_{lc} & [0]_{2p \times (N+1-4p)} \\ [0]_{(N+1-4p) \times (p+1)} & D^-_c & [0]_{(N+1-4p) \times (p-1)} \\ [0]_{2p \times (N+1-4p)} & D^-_{rc} & D^-_r \end{pmatrix}, \qquad (12)$$

where the stencil $S_i : j = i - p - 1, ..., i + p - 1$ is used for $D^-_c$.

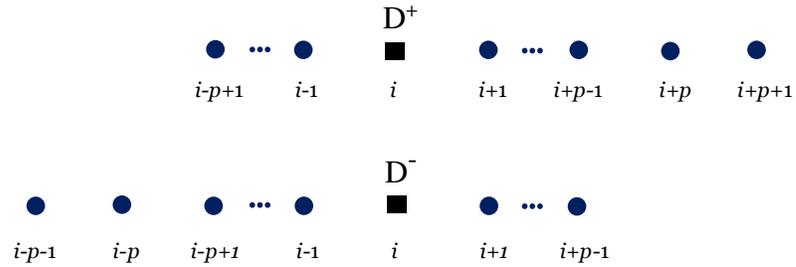

Fig. 2 The operator stencils in IGNs

Consider the matrices $D_c^+$, $D_c^-$, see. Fig. 2. The coefficients $\beta_{-p+1},...,\beta_{p+1}$ for the elements of $D_c^+$ are evaluated by the standard way of Taylor series expansion of a smooth function $f(x)$ in the neighborhood of $i$ node for the shifted stencil $S_i : j = i - p+1,...,i+p+1$. Namely, $2p$ equations are written out to provide $2p$-order approximation:

$$f_{j+i} = f_i + f'_i jh + b_1(jh)^2 + ... + b_{2p-1}(jh)^{2p}, \quad j = -p+1,...,p+1, \quad j \neq 0.$$

Thus, we obtain the corresponding values of $\beta_j$ as coefficients at $f_{j+i}$. Similarly, for the entries of $D_c^-$.

Now consider the matrices $D_{lc}^+$, $D_{lc}^-$. The coefficients $\{\alpha_{j,k}\}$ for the entries of $D_{lc}^+$ are evaluated from (5):

$$D^+ = -H^{-1}(D^-)^T H + H^{-1}Q. \tag{13}$$

The matrix $D_{lc}^+$ corresponds to the submatrix of the right hand side in (13) with the numbers of rows $p,...,2p$ and columns $2p+1,...,3p+1$, i.e., elements of the matrix $D_c^-$ that are known. Therefore, the matrix $D_{lc}^+$ is uniquely determined by elements of $H$, i.e., the functional relationship

$$\alpha_{j,k} = \alpha_{j,k}(\mu_1,...\mu_{2p}). \tag{14}$$

is known. Similarly, the matrix elements of $D_{lc}^-$ are parameterized in terms of $(\mu_1,...\mu_{2p})$ of $H$.

Now consider the matrices $D_l^+$, $D_l^-$ and $H$. We require that the approximation of the first derivative in the left NBNs was of $p_a \geq p$ order. To do this, we write the appropriate conditions for monomials $x^n$, $n = 0,1,...$, in the points (1) for the corresponding $2p$ points stencils:

$$\begin{cases} D_l^+ w_n = n w_{n-1} \\ D_l^- w_n = n w_{n-1} \end{cases}, \quad n = 0,...,p_a, \tag{15}$$

where the vector $w_n$ is defined by $[w_n]_{2p \times 1} = (x_0^n,...,x_{2p-1}^n)^T$.

The equation (5) allows us to express the elements of $D_l^-$ in terms of $D_l^+$ and parameters $\mu_1,...,\mu_{2p}$ and $h_1, h_2,..., h_{p_g-1}$. Therefore, considering (15) as a system for entries of $D_l^+$ we get $2 \times 2p \times (p_a + 1)$ equations, whereas the number of unknowns is just $2p \times 2p = 4p^2$. I.e., we obtain the overdetermined system. To ensure its compatibility it is necessary to find appropriate parameters $\mu_1,...,\mu_{2p}$ and $h_1, h_2, ..., h_{p_g-1}$. We use the following approach. Let us set $p_a = p$. The overdetermined system (15) is analytically transformed to the upper triangle matrix in the first $4p^2$ rows. Denote this matrix by $UT(d \mid \mu, h)$. The remaining $4p$ rows form the equations with zero left hand sides and some right hand sides linearly dependent on $\mu_1, ..., \mu_{2p}$ and nonlinearly dependent on $h_1, h_2, ..., h_{p_g-1}$. Denote the system of these equations by $DS(\mu \mid h)$. One can solve $DS(\mu \mid h)$ with respect to $\mu_1, ..., \mu_{2p}$ after taking some $h_1, h_2, ..., h_{p_g-1}$. Our numerical experiments have identified the following features of this system for $p = 2, ..., 6$, ($p_a = p$):

**Observation 1.** Despite the double overdetermination of $DS(\mu \mid h)$ with respect to $\mu_1, ..., \mu_{2p}$ ($4p$ equations for $2p$ unknowns), the rank of the system is $2p$, i.e., the system always has a solution.

**Observation 2.** For the uniform grid $h_1 = ... = h_{p_g-1} = h$ the criterion $\mu_i > 0$, see (6), is satisfied for $p = 2, 3, 4$ only.

**Observation 3.** There are sets $h_1, h_2, ..., h_{p_g-1}$ for grids with shifted NBNs when $\mu_i > 0$ also at $p = 5, 6$.

**Observation 4.** After substituting the determined $\mu_1, ..., \mu_{2p}$ from solution of $DS(\mu \mid h)$ in the matrix $UT(d \mid \mu, h)$ its rank becomes equal to $4p^2 - (p-1)^2$, i.e., the coefficients $\{d_{i,j}\}$ linearly depend on another $(p-1)^2$ free parameters $c_1, c_2, ..., c_{(p-1)^2}$.

It should be noted that any attempts to consider $p_a > p$ lead to a violation of observation 1: the system $DS(\mu \mid h)$ becomes insoluble.

Summing up, we obtain the following experimental numerical results: for $p = 2, ..., 6$ and $p_a = p$ at some $h_1, h_2, ..., h_{p_g-1}$, $p_g \leq p$, there exist positive sets $\mu_1, ..., \mu_{2p}$ (i.e., the matrix $H$) and coefficients of matrix $D_l^+$, linearly dependent on $(p-1)^2$ free parameters:

$$d_{i,j} = d_{i,j}\left(c_1, ..., c_{(p-1)^2}\right), \quad i, j \in [1, ..., 2p]. \tag{16}$$

Recall, that coefficients of the matrix $D_l^-$ are evaluated from $D_l^+$ and $H$ using (5).

The similar procedure is applied to the bottom third of the matrices $D^+, D^-$, i.e., for the matrices $D_r^+$, $D_r^-$, $D_{rc}^+$, $D_{rc}^-$.

## 1.4 Search for optimal parameters of grid and operators $D^+, D^-, H$

Let us will choose free parameters $\vec{h} \equiv (h_1, h_2, ..., h_{p_g-1})$ and $\vec{c} \equiv (c_1, ..., c_{(p-1)^2})$ to fulfill the conditions of the second part of **P2** and **P3**.

To do this, we formulate the problem of minimizing the objective functional consisting of two components responsible for the **P1** and **P2**:

$$\mathrm{E}(\vec{h},\vec{c}) \equiv \sum_{n=p+1}^{2p} \frac{\left\| D^- D^+ T_n - T_n'' \right\|_H^2}{\left\| T_n \right\|_H^2} + C\theta\left(\frac{\lambda_{full}}{\lambda_{int}} - \kappa\right) \to \min_{\vec{h},\vec{c}}. \qquad (17)$$

Here $T_n = (T_n(x_0), T_n(x_1), ..., T_n(x_N))^\mathrm{T}$; $T_n'' = (T_n''(x_0), T_n''(x_1), ..., T_n''(x_N))^\mathrm{T}$; $T_n(x)$ are the Chebyshev first kind polynomials; $C > 0$, $\kappa > 1$ are some constants; $\theta(x)$ is the Heaviside function; $(\lambda_{full})^2$ and $(\lambda_{int})^2$ are the maximal eigenfunctions of the problem $(D^+)^\mathrm{T} H D^+ \bar{u} = \lambda^2 H \bar{u}$ on the whole grid and on IGNs, respectively (for the periodic functions in the latter case).

It should be noted that due to the condition **P1** the equality $D^- D^+ T_n - T_n'' = 0$ is valid at IGNs up to $n=2p$. Eigenvalues $\lambda_{int}$ are computed in advance for periodical problem on uniform grid and do not depend on $\vec{h}, \vec{c}$ during optimization problem.

Since $\mathrm{E}(\vec{h}, \vec{c})$ depends on its variables in a complex way, it was used the following algorithm to find them (maybe not the best, but it led to success).

We introduce the auxiliary functional

$$\mathrm{E}_{aux}(\vec{h}, \vec{c}) \equiv \sum_{n=p+1}^{2p} \left\| D^+ T_n - T_n' \right\|_H^2. \qquad (18)$$

We take a set $\vec{h}$, and consider the minimization problem of $\mathrm{E}_{aux}(\vec{h}, \vec{c})$ with respect to $\vec{c}$. Its solution obviously reduces to solving a system of linear equations with respect to $\vec{c}$ because of the linear dependency of $D^+$ on $\vec{c}$ as previously indicated. Let us denote this operation by $C_{aux}$, i.e., $\vec{c} = C_{aux}(\vec{h})$. Now, by taking some initial guess for $\vec{h}$, we find the solution

$$\vec{h}_0 = \arg\min_{\vec{h}} \ \mathrm{E}\left(\vec{h}, \ \vec{c} = C_{aux}(\vec{h})\right). \tag{19}$$

using the Nelder-Mead simplex method (see., e.g., [2]).

Afterwards we sequentially solve the following problems by the same simplex method:

$$\{\vec{c}_i: \ \mathrm{E}\left(\vec{h}_{i-1},\vec{c}\right) \to \min_{\vec{c}}; \quad \vec{h}_i: \ \mathrm{E}\left(\vec{h},\vec{c}_i\right) \to \min_{\vec{h}}\}, \quad i=1,2,.... \tag{20}$$

Typically, the process (20) finishes at $i=2$ with satisfactory enough results.

## 2 Results

Here we report the results of testing the developed difference schemes on the one-dimensional problem for the wave equation.

### 2.1 Test problem

To identify the accuracy of the proposed schemes, we consider the initial boundary value problem for the wave equation on the interval $x \in [-0.5, 0.5]$ with homogeneous Neumann boundary condition and with a Gaussian profile as the initial condition, see Figure 3:

$$\begin{cases} \dfrac{\partial^2 u}{\partial t^2} - \dfrac{\partial^2 u}{\partial x^2} = 0 \\ u\big|_{t=0} = e^{-0.5(x/0.05)^2}, \quad u'_t\big|_{t=0} = u'_x\big|_{x=-0.5} = u'_x\big|_{x=0.5} = 0 \end{cases}.$$

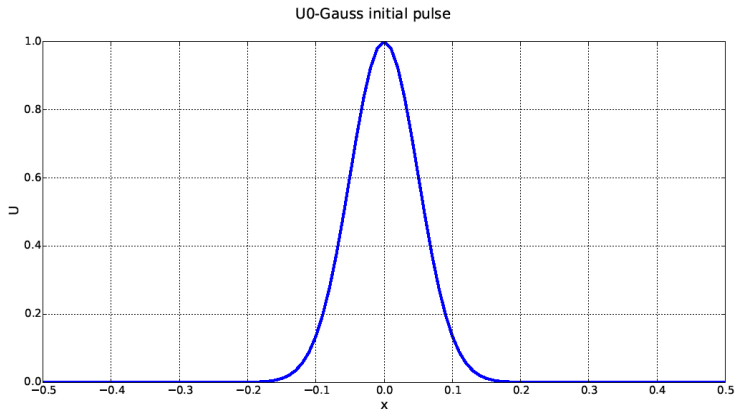

**Fig. 3 The initial data**

This task permits to evaluate the accuracy of the schemes both before the impact of boundaries, i.e., at $t < 0.3$, and after the instance when the impulse comes to the ends of the segment.

We use grids with the nodes number $N+1 = 101, 111, 121, 131, 151, 171, 201, 231, 261, 301$.

## 2.2 Accuracy of the schemes

Here we present some characteristics of several finite-difference schemes based on the approach suggested above.

| $2p$ | $K$ | $P_{num}$ | $\lambda_{int}/\lambda_{full}$ | $2h/\lambda_{int}$ |
|---|---|---|---|---|
| 4 | 0 | 4.3 | 0.46 | 0.75 |
| 4 | 1 | 3.4 | 0.54 | 0.75 |
| 6 | 0 | 4.4 | 0.30 | 0.91 |
| 6 | 1 | 6.1 | 0.29 | 0.91 |
| 6 | 2 | 4.9 | 0.26 | 0.91 |
| 8 | 0 | 4.2 | 0.23 | 0.95 |
| 8 | 1 | 5.1 | 0.24 | 0.95 |
| 8 | 2 | 6.4 | 0.13 | 0.95 |
| 8 | 3 | 4.2 | 0.16 | 0.95 |
| 10 | 1 | 5.7 | 0.22 | 0.96 |
| 10 | 2 | 8.2 | 0.12 | 0.96 |
| 12 | 1 | 5.6 | 0.19 | 0.95 |
| 12 | 2 | 9.0 | 0.13 | 0.95 |

**Tab. 1 Main information of the generated schemes**

Table. 1 contains information about

$2p$, the order of the scheme

$K = p_g - 1$, the number of shifted NBNs (0 means the uniform grid)

$P_{num}$, the observed numerical order of accuracy at $t = 0.5$ (i.e., when it is significantly influenced by boundaries)

$\lambda_{int} / \lambda_{full}$, the coefficient of reduction of the Courant number because of the influence of NBNs, and

$2h / \lambda_{int}$, the theoretical Courant number in IGNs for the explicit second order scheme in time.

Column with $\lambda_{int} / \lambda_{full}$ refers to matrices on the grid with $N + 1 = 101$. The coefficients of some of the constructed schemes are given in Annex A.

Figure 4 shows graphs of the grid convergence to an exact solution for the schemes with $2p = 8$ at $t = 0.2$ (no influence of the boundaries, therefore, $P_{num} = 7.6$ in all cases) and at $t = 0.5$.

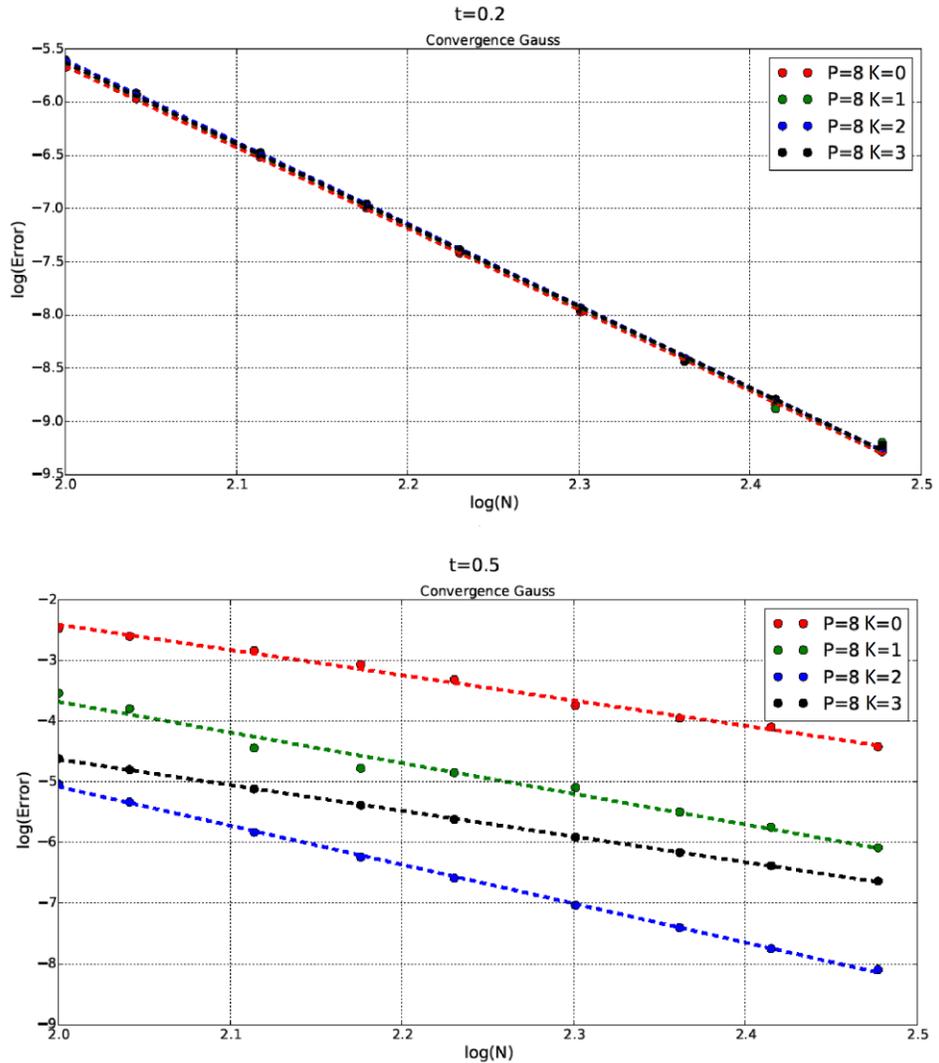

**Fig. 4 The grid convergence for schemes with *2p=8* at t=0.2 (top) and t=0.5 (bottom)**

The next series of tests allows comparing the accuracy characteristics of the proposed schemes. All calculations use the grid with $N = 101$. Consider the case $2p = 8$. Figure 5 shows plots of numerical solution errors depending on time. We see, for example, that the scheme on the equidistant grid has approximately 100 times greater error than the scheme on the grid with two shifted NBNs ($K = 2$).

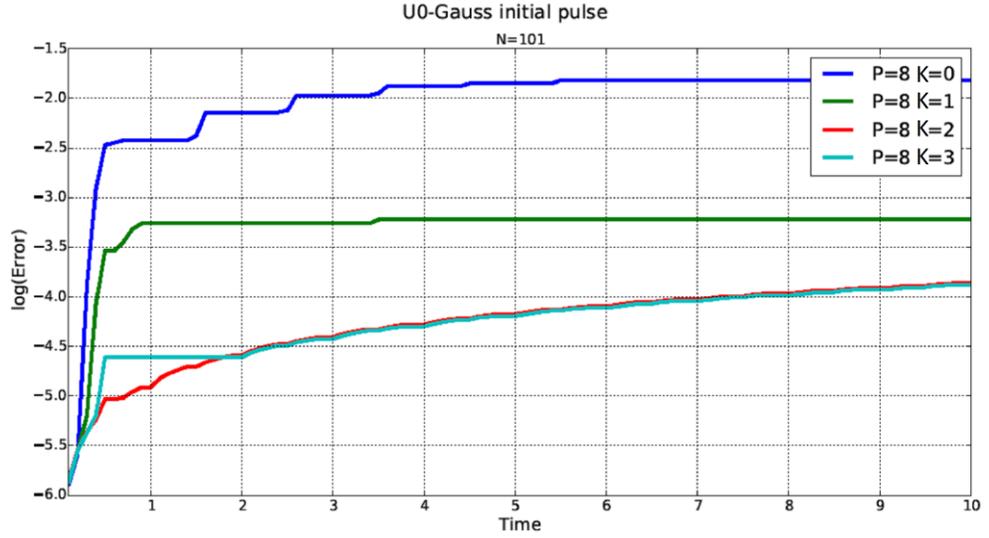

**Fig. 5 Solution error vs time in C-norm for schemes with *2p=8***

Now compare the accuracy of schemes of different orders: $2p = 4,..., 12$. Fig. 6 shows graphs of the numerical solution errors for values of $K$ that provide the maximum accuracy. One can see that the best accuracy was obtained for the scheme with $2p = 10$, $K = 2$. A low accuracy of the scheme with $2p = 12$ is explained by unsuccessful optimization (17) of its coefficients due to much more computationally complex tasks compared to $2p = 10$. Further work is needed to improve the algorithm of solving (17).

We note that according to Table 1 the Courant number of the difference scheme with, e.g., $2p = 10$, $K = 2$, is reduced by factor about *eight* at NBNs ($0.96/0.12 \approx 8$). This may seem a rather serious shortcoming of the proposed schemes. However, such many-fold reduction in Courant number for the one-side stencil schemes at the end points of the segment is quite natural thing. Moreover, the Courant number of the conventional SBP scheme on the equidistant grid with $2p = 8$, $K = 0$ is reduced by factor about *four* ($0.95/0.23 \approx 4$). However, the accuracy of the new scheme is about 500 times as better, see Figs. 5 and 6. Notice that we found another variant of scheme with the same values of order $2p=10$ and $K=2$, which has just *three-fold* reduction in Courant number at NBNs, but its accuracy is about 50 times as better as that of the abovementioned conventional SBP scheme. We highlight that there is a choice of priority criteria either for the accuracy or for the Courant number.

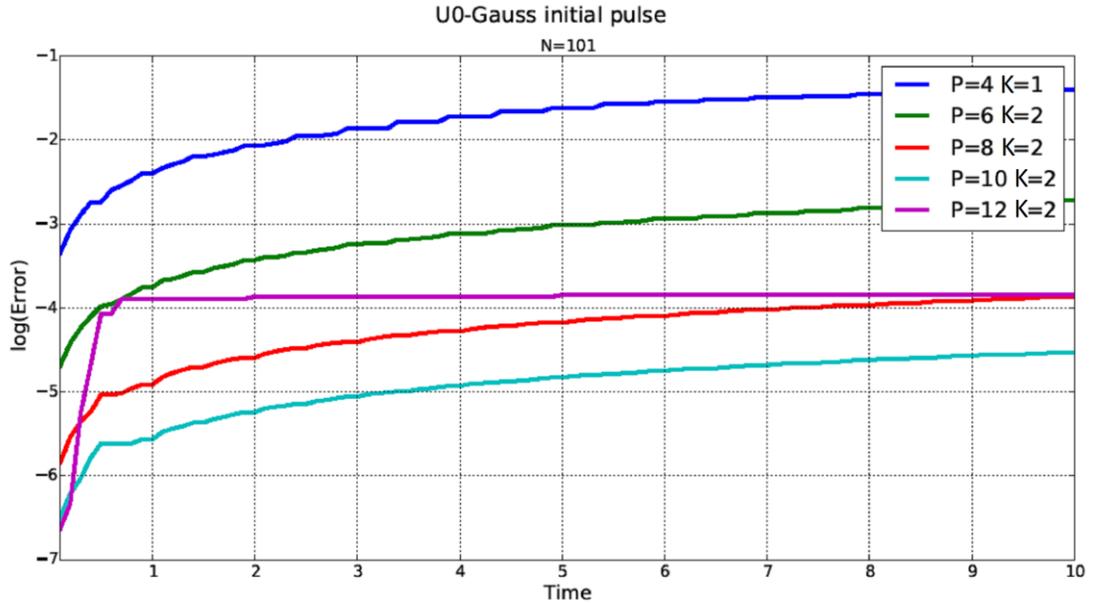

**Fig. 6 Solution error vs time in C-norm for best schemes with *2p=4, ..., 12***

## Conclusion

High-order finite-difference operators for the first and second derivatives on the segment are proposed. The corresponding algorithm for generating their coefficients is developed.

The schemes are based on:

- SBP approach with forward and backward differences
- use of shifted near boundary nodes with respect to the equidistant grid
- sequential application of the first derivative operators to obtain finite-difference operators of the second derivative

The proposed difference operators have the following properties:

- high-order approximation (up to 12th)
- symmetric matrices for the second derivative operators
- damping of saw-tooth spurious oscillations in solutions

Essentially higher accuracy is observed for the new schemes comparing with the schemes on the equidistant grid while solving test problems for the wave equation. It was shown that the natural order decay due to SBP boundary closure is less for the new schemes: e.g. for the $8^{th}$ order scheme with two shifted nodes ($2p=8, K=2$) it reduces to 6.4 instead of 4.2 for ($2p=8, K=0$), the Neumann problem.

New schemes with $2p=8, K=2$ and $2p=10, K=2$ demonstrate best accuracy among other constructed ones.

An important feature of the new schemes for the wave equation is that the Courant number is not less than half of that for the schemes on the equidistant grid, but the accuracy is more than 100 times as better in considered tests. To our knowledge this is the first examples of SBP finite difference operators of $10^{th}$ and $12^{th}$ order with a diagonal normalization matrix $H>0$, which are presented in literature.

# Appendix. Difference operators for grids with shifted NBNs (K=1,2,3)

## 1. Grid nodes

Here the spacing $h = 1$ in IGNs

| $2p$ | $K=1$ | $K=2$ | $K=3$ |
|---|---|---|---|
| 4 | $h_1 = 0.64701892044823239$<br>$h_2 = h = 1$ | - | - |
| 6 | $h_1 = 0.55959440808516225$<br>$h_2 = h = 1$ | $h_1 = 0.52989554067209088$<br>$h_2 = 0.9577049256058392$<br>$h_3 = h = 1$ | - |
| 8 | $h_1 = 0.53057599940567612$<br>$h_2 = h = 1$ | $h_1 = 0.39203322551059488$<br>$h_2 = 0.81423930361885499$<br>$h_3 = h = 1$ | $h_1 = 0.43979786646687147$<br>$h_2 = 0.90985090947051206$<br>$h_3 = 1.0771428495647428$<br>$h_4 = h = 1$ |
| 10 | $h_1 = 0.50900297608285072$<br>$h_2 = h = 1$ | $h_1 = 0.37366515483267776$<br>$h_2 = 0.79308655639992476$<br>$h_3 = h = 1$ | - |
| 12 | $h_1 = 0.48125000596046169$<br>$h_2 = h = 1$ | $h_1 = 0.38823311074361344$<br>$h_2 = 0.81640993512856175$<br>$h_3 = h = 1$ | - |

## 2. Matrix of operator $H$ and free parameters

**1.** $2p = 4$

1.1. $K=1$

$\mu_1 = 0.186109276322411116, \mu_2 = 0.975448598874482986, \mu_3 = 0.976489275826791681$
$\mu_4 = 1.008971769424537701, \mu_5 = 1$

**2.** $2p = 6$

2.1. $K=1$

$\mu_1 = 0.162227980272819955, \mu_2 = 0.873555067807182617, \mu_3 = 1.031381558634351325$
$\mu_4 = 0.991130867107816504, \mu_5 = 1.001296776753106244, \mu_6 = 1.000002157509889189$
$\mu_7 = 1$

2.2. *K*=2

$\mu_1 = 0.153545834255111785, \mu_2 = 0.827868630728788024, \mu_3 = 1.007836990306931968$
$\mu_4 = 0.998846303560314341, \mu_5 = 0.999229946106053313, \mu_6 = 1.000272761320736059$
$\mu_7 = 1$

3. $2p = 8$

3.1. *K*=1

$\mu_1 = 0.294839655769715270, \mu_2 = 1.526077766754849963, \mu_3 = 0.256381448412698443$
$\mu_4 = 1.799899415784832479, \mu_5 = 0.410922343474426854, \mu_6 = 1.279556051587301679$
$\mu_7 = 0.922938436948853691, \mu_8 = 1.009384881267321621, \mu_9 = 1$

3.2. *K*=2

$\mu_1 = 0.110338815724131761, \mu_2 = 0.635841857271623012, \mu_3 = 0.950804714528380446$
$\mu_4 = 1.013133760425796615, \mu_5 = 0.994604889106195045, \mu_6 = 1.001983074003966800$
$\mu_7 = 0.999507111548435856, \mu_8 = 1.000058306520917872, \mu_9 = 1$

3.3. *K*=3

$\mu_1 = 0.123920978343533647, \mu_2 = 0.712830551002903490, \mu_3 = 1.054602637461168113$
$\mu_4 = 1.046653662344426694, \mu_5 = 0.985006468579236016, \mu_6 = 1.004695892473460361$
$\mu_7 = 0.998971519552118048, \mu_8 = 1.000109915745287070, \mu_9 = 1$

4. $2p = 10$

4.1. *K*=1

$\mu_1 = 0.144437776901122500, \mu_2 = 0.817330306276763396, \mu_3 = 1.085198183908459013$
$\mu_4 = 0.929611527349792244, \mu_5 = 1.055478528306474040, \mu_6 = 0.964711409696087818$
$\mu_7 = 1.016688169724728752, \mu_8 = 0.994562982288813791, \mu_9 = 1.001083360908970654$
$\mu_{10} = 0.999900730721632103, \mu_{11} = 1$

4.2. *K*=2

$\mu_1 = 0.104654633738292063, \mu_2 = 0.609748535967006844, \mu_3 = 0.940178174636667863$
$\mu_4 = 1.018013035012260925, \mu_5 = 0.991135498769543100, \mu_6 = 1.004360329350681758$
$\mu_7 = 0.998221653537359699, \mu_8 = 1.000531902166808873, \mu_9 = 0.999898886800964282$
$\mu_{10} = 1.000009061253016585, \mu_{11} = 1$

5. $2p = 12$

    5.1. *K=2*

$\mu_1 = 0.108740366453879106, \mu_2 = 0.632865778556215175, \mu_3 = 0.956449007732071865$
$\mu_4 = 1.007176273126186183, \mu_5 = 1.001241188722324038, \mu_6 = 0.995791403173844181$
$\mu_7 = 1.004146733839414329, \mu_8 = 0.997356904816331769, \mu_9 = 1.001160562212721983$
$\mu_{10} = 0.999659249972080222, \mu_{11} = 1.000060503196924522, \mu_{12} = 0.999995074070187173$
$\mu_{13} = 1$

The coefficients of $D^+, D^-$ are computed by using the above sets of $\vec{h}$, $\vec{\mu}$.